\documentclass[12pt]{amsart}

\usepackage{amsthm}
\usepackage{amsmath}
\usepackage{amssymb}
\usepackage{graphicx}
\usepackage{enumerate}
\allowdisplaybreaks
\numberwithin{equation}{section}
\numberwithin{figure}{section}

\theoremstyle{plain}
\newtheorem{theorem}{\sffamily Theorem}[section]
\newtheorem{proposition}[theorem]{\sffamily Proposition}
\newtheorem{lemma}[theorem]{\sffamily Lemma}
\newtheorem{corollary}[theorem]{\sffamily Corollary}
\newtheorem{example}[theorem]{\sffamily Example}
\newtheorem{remark}[theorem]{\sffamily Remark}
\newtheorem{definition}[theorem]{\sffamily Definition}
\newtheorem{conjecture}{\sffamily Conjecture}

\def\BET{\begin{theorem}}
\def\ENT{\end{theorem}}
\def\BEP{\begin{proposition}}
\def\ENP{\end{proposition}}
\def\BEL{\begin{lemma}}
\def\ENL{\end{lemma}}
\def\BEC{\begin{corollary}}
\def\ENC{\end{corollary}}
\def\BEE{\begin{example}\rm}
\def\ENE{\end{example}}
\def\BER{\begin{remark} \rm}
\def\ENR{\end{remark}}
\def\BED{\begin{definition} \rm}
\def\END{\end{definition}}
\def\BECJ{\begin{conjecture}}
\def\ENCJ{\end{conjecture}}

\def\bea{\begin{eqnarray}}
\def\eea{\end{eqnarray}}

\def\beas{\begin{eqnarray*}}
\def\eeas{\end{eqnarray*}}

\def\beq{\begin{equation}}
\def\eeq{\end{equation}}

\def\beal{\begin{align*}}

\def\eeal{ \end{align*} }

\def\row{ \nonumber \\ & & }
\def\roweq{\nonumber \\ &=& }
\def\rowleq{\nonumber \\  & \leq & }
\def\rowgeq{\nonumber \\ & \geq & }

\def\rowpl{\nonumber \\ & \ \ &+  }

\def\bbC{{\mathbb C}}
\def\bbD{{\mathbb D}}

\def\bbN{{\mathbb N}}

\begin{document}

\title[Boundedness of Toeplitz operators]{On the boundedness of Toeplitz operators with radial symbols over weighted sup-norm spaces of holomorphic functions}

\author{Jos\'e Bonet}
\author{Wolfgang Lusky}
\author{Jari Taskinen}

\thanks{Acknowledgement:  
The research of Bonet was partially supported by the projects MTM2016-76647-P and GV Prometeo/2017/102.} 

\begin{abstract}
We prove sufficient conditions for the boundedness and compactness of
Toeplitz operators $T_a$ in weighted sup-normed Banach spaces $H_v^\infty$
of holomorphic functions defined on the open unit disc $\bbD$ of the complex
plane; both the weights $v$  and symbols $a$ are assumed to be radial
functions on $\bbD$. In an earlier work by the authors it 
was shown that there exists a bounded, harmonic (thus non-radial)
symbol $a$ such that $T_a$ is not bounded in any space $H_v^\infty$ with
an admissible weight $v$. Here, we show that a mild additional assumption
on the logarithmic decay rate of a radial symbol $a$  at the boundary of $\bbD $ 
guarantees the boundedness  of $T_a$.  

The sufficient conditions for the boundedness and compactness of $T_a$, 
in a number of variations, are derived from the general, 
abstract necessary and  sufficient condition recently found by the authors. 
The results apply for a large class of weights satisfying the so called condition
$(B)$, which includes in addition to standard weight classes also many
rapidly decreasing weights.
\end{abstract}

\maketitle

\section{Introduction and main results.}
\label{sec1}
In the article \cite{toe} we studied Toeplitz operators $T_a$ with radial 
symbols $a$  on the analytic function spaces $H_v^\infty$ on the unit 
disc $\bbD \subset \bbC$, endowed with weighted sup-norms for a large class 
of  radial weights $v$  satisfying the so called condition $(B)$; this 
excludes the unweighted or constant weight case. In particular, 
in Theorem 3.6 of the citation (repeated in this paper in Theorem 2.1) 
we obtained a general sufficient and necessary condition for the  
boundedness and compactness of $T_a: H_v^\infty \to H_v^\infty$. 
Also, we observed that the boundedness of a non-radial symbol does not necessarily  
imply the boundedness of the Toeplitz operator. 
In fact, Theorem 2.3 of \cite{toe} contains an example of a 
bounded harmonic symbol $a$  such that $T_a : H_v^\infty \to H_v^\infty$ is 
not bounded for any weight $v $ under consideration. 

The criterion for the boundedness of the Toeplitz operator in Theorem 3.6 of  
\cite{toe} is quite abstract, and it may
not be easy to verify it for concrete weights and symbols. Some examples
were presented in the citation under quite special assumptions
either on the symbol or on the weight.
Here, our aim is to use 
Theorem 3.6 of \cite{toe} to prove concrete sufficient conditions for the 
boundedness and  compactness  of $T_a : H_v^\infty \to H_v^\infty$.  These 
conditions are much more general than in the examples of the citation,  and the 
sufficient conditions for the symbol are easy to 
formulate and control. In all of our results we assume that the 
weight $v$ satisfies condition $(B)$ of \cite{lu2}, see also Definition 
\ref{def1.1} below, and a mild technical condition \eqref{1.1}. These 
assumptions hold for example for  the important classes of standard, normal 
and exponential weights  (Proposition \ref{prop1.3}). Then, in the first main
result,  Theorem \ref{th1.2}, we show that for the boundedness of 
$T_a : H_v^\infty \to H_v^\infty$ it suffices that the symbol $a$
is differentiable near the unit circle and $\limsup a'$ 
or $- \limsup a'$ 
is bounded from above, and $a \to 0$ at a slow, logarithmic speed as 
$r \to 1$. In Theorem \ref{cor1.4} the decay requirements for $a$ are replaced 
by decay conditions on $a'$. In the case of normal weights the smoothness requirements of the symbol can be relaxed, see Theorem \ref{th1.5}. 
Finally, in Theorem \ref{th1.6} we find a stronger decay 
condition for $a$ which guarantees the boundedness of $T_a$ in the case of 
exponential weights without a smoothness assumption on $a$.

All of these theorems also contain the analogous statements on the 
compactness of the Toeplitz operator. The proofs of Theorems
\ref{th1.2} and \ref{th1.5} will be presented in Section \ref{sec3}
and that of Theorem \ref{th1.6} in Section \ref{sec4}.

We refer to the papers   \cite{E}, \cite{G},  \cite{G1},  \cite{Lu},
\cite{Lu2},  \cite{LT},  \cite{M}, \cite{PTV1},  \cite{S},  \cite{SZ}, \cite{Su}, \cite{tasvir}, \cite{V}, \cite{V1}, \cite{Z1}, \cite{zh}, \cite{Zo} for classical and recent results on the boundedness and  compactness of Toeplitz operators on Bergman spaces

Let us turn to the exact definitions and  formulation of the main results. 
By  a weight $v$ on the unit disc $ \mathbb{D}$ we mean  a continuous  
function  with $v(z)= v(|z|)$ for all $z \in \mathbb{D}$, $\lim_{|z| \rightarrow 1}v(z) =0$ and $v(r) \geq v(s)$ if
$1 > s > r >0$. Put
\beas 
& & H^{\infty}_v= \big\{ h: \mathbb{D} \rightarrow \mathbb{C} : h \mbox{ holomorphic}, \ \Vert h \Vert_v := \sup_{z \in  \mathbb{D}} |h(z)|v(|z|) < \infty \big\}, 
\row
L^{\infty}_v= \big\{ h: \mathbb{D} \rightarrow \mathbb{C} : h \mbox{
measurable }, \ \Vert h \Vert_v := {\rm ess}\sup_{z \in  \mathbb{D}} |h(z)|v(|z|) < \infty \big\}.  
\eeas
Let $\mu$ be the Lebesgue area measure on $ \mathbb{D}$ endowed with $v$ as density, i.e. $d\mu(r e^{i \varphi}) = v(r)rdr d\varphi$ and denote the weighted
$L^p$- and Bergman spaces by
 \[ 
L^p_v = \Big\{ g:  \mathbb{D} \rightarrow  \mathbb{C} : \Vert g \Vert_{p,v}^p := 
\int_{ \mathbb{D}} |f|^p d\mu < \infty \Big\}  \mbox{ and }  A_v^p = \{ h \in L^p_v : h \mbox{ holomorphic} \},
 \] 
where $1 \leq p < \infty $. In the unweighted case $v$ is omitted in the 
notation.
 
Now let $a  \in L^1$. We define the Toeplitz operator $T_a$ with 
symbol $a$ on  $H^{\infty}_v$  by $T_a h = P_v(a \cdot h)$ for $h \in  
H^{\infty}_v $,   where $P_v :L_v^2 \rightarrow H_v^2$ is the orthogonal 
projection. Then $T_a h$ is a holomorphic function, at least if 
$a \cdot h \in L_v^2$. The definition of the Toeplitz operator in the 
present setting is discussed in detail in  Section 1 of  
\cite{toe} and we do not wish to repeat the details here. However, we 
emphasize that even if  $T_a h$ is a well a well defined analytic 
function, is not necessarily an element of 
$H^{\infty}_v$ and $T_a$ need not be a bounded operator.
  
In the following we consider radial symbols $a \in L^1$, i.e. 
functions with  $a(z) = a(|z|)$ for almost all $z \in\mathbb{D}$.
As for general notation, 
$\bbN  = \{ 1,2,3, \ldots \}$, $\bbN_0 = \bbN \cup \{ 0 \}$, and $c,C,C'$ 
denote generic positive constants, the exact value of which may change from 
place to place, but does not depend on the variables, indices or 
functions in the given expressions, unless otherwise indicated. By 
$1_A$ we denote the characteristic function of a set $A$, i.e. 
a function which equals 1 on $A$ and $0$ outside $A$; the domain of $1_A$ 
will be clear from the context. For other general terminology and definitions,
see \cite{hkz} and \cite{zh}.

We will need the following definition. Let $v$ be a weight on $\mathbb{D}$. 
Consider $m >0$ and let $r_m$ be a point where the function
$r^mv(r)$ attains its absolute maximum on $[0,1]$. It is easy to see that $r_n \geq r_m$ if $n \geq m$ and $\lim_{m \rightarrow \infty} r_m = 1$.

\noindent
\BED \label{def1.1}
(i) The weight $v$ satisfies  condition $(B)$, if
\begin{eqnarray*}
	& & \forall b_1 >1 \  \exists b_2 >1 \  \exists c>0 \
	\forall m, n  > 0  \\
	& & \left( \frac{r_m}{r_n} \right)^m \frac{v(r_m)}{v(r_n)}
	\leq b_1
	\ \ {\rm and} \ \  m , n , |m-n| \geq c \ \
	\Rightarrow
	\left( \frac{r_n}{r_m} \right)^n \frac{v(r_n)}{v(r_m)} \leq b_2,
\end{eqnarray*}

\noindent (ii)  $v$ is called {\em normal} if
\[
\sup_{n\in \mathbb{N}} \frac{v(1-2^{-n})}{v(1-2^{-n-1})} < \infty \ \ \
\mbox{ and } \ \ \ \inf_k \limsup_{n \rightarrow \infty}\frac{v(1-2^{-n-
k})}{v(1-2^{-n})} < 1,
\]

\noindent  (iii) $v$ is called an {\em exponential weight of type}
$(\alpha,\beta)$ for some  constants $\alpha > 0$ and $\beta > 0$ if
$v(r) = \exp(- \alpha/(1-r)^\beta)$ for all $r$.
\END

Note  that the numbers $m$ and $n$ in (i) need not be integers.
For example all normal weights as well as all  exponential weights
satisfy $(B)$  (see \cite{lu2}). Standard weights $ (1-r)^\alpha$
and $ (1-r^2)^\alpha$ are normal for all $\alpha > 0$, but no  exponential 
weight is normal: the  first condition in (ii) is not satisfied. Neither is
the weight  $v(r) = 1/(1-  \log(1-r))$ is normal, since it decays too slowly 
to 0 in order to satisfy the second condition in (i).   

We show

\BET
\label{th1.2}
Let $v$ satisfy $(B)$ and assume that there is some $\epsilon > 0$ with
\bea \label{1.1} 
\sup_{n=1,2,\ldots}\frac{\int_0^1r^{n-n^{\epsilon}}v(r) dr}{\int_0^1r^nv(r)dr} < \infty.
\eea
Let $a \in L^1$ be real valued and radial  such that the restriction 
$a|_{[ \delta,1[} $ is  differentiable for some $\delta \in ]0,1[$ with 
\bea
\limsup_{r \rightarrow 1} a'(r) < \infty \ \ \ \ \ \ \mbox{ or } \ \ \ \ \ \ \liminf_{r \rightarrow 1}a'(r) > - \infty.  \label{1.1a}
\eea 
If
\bea \label{1.2}
\limsup_{r \rightarrow 1} |a(r) \log(1-r)| < \infty
\eea
then $T_a$ is a bounded operator $  H^{\infty}_v \rightarrow  H^{\infty}_v$.

If
\bea
   \limsup_{r \rightarrow 1}  |a(r) \log(1-r) | = 0
\label{1.3}
\eea
then $T_a$ is compact on $H_v^{\infty}$.
\ENT

Of course, Theorem 1.2 can be applied to complex valued symbols $a$ as well. Here 
Re $a$ and Im $a$ have to satisfy the assumptions of the theorem.

We prove Theorem \ref{th1.2} in Section \ref{sec3}. Condition  $(B)$ and 
\eqref{1.1} are satisfied for many weights, in particular we have 

\BEP \label{prop1.3}
All normal and exponential weights (see Definition \ref{def1.1}) satisfy $(B)$ and  
condition \eqref{1.1}. 
\ENP

Indeed, it was proven in \cite{lu2} that normal and exponential
weights satisfy $(B)$. Condition \eqref{1.1} with $\epsilon = 1/2$ follows for 
normal weights from Lemma 4.5. of [3]. The remaining claim  in Proposition 
\ref{prop1.3} about condition \eqref{1.1} for exponential weights will be 
proved in  Section \ref{sec4}; see the remark after the proof of Corollary 
\ref{cor4.2}. 

\bigskip

{\bf Examples.}
Assume that $v$ is a weight satisfying $(B)$.
The symbol  $a(r) = 1/(1-\log(1-r))$ satisfies the second condition 
\eqref{1.1a} and, of course, \eqref{1.2} so that  $T_a: H_{v}^{\infty} \to  
H_{v}^{\infty}$ is bounded. The same is true for $a(r)= (1-
r)^{\delta}$ with  any $\delta > 0$, and this symbol even satisfies 
\eqref{1.3} so that $T_a$ is compact on  $H_{v}^{\infty}$.  (1.3). Moreover, 
\[ a(r) = \left\{ \begin{array}{rc}
\log 2 \, , & \mbox{ if } \ 0 \leq r \leq 1/2, \\
-\log r \, , & \mbox{ if }  \ 1/2 < r < 1,
\end{array} \right. \]
satisfies \eqref{1.1a}, \eqref{1.3}  as well.

Next we present a reformulation of Theorem \ref{th1.2}.
\BET \label{cor1.4}
Let $v$ satisfy $(B)$ and \eqref{1.1}. Let  $a \in L^1$ be a radial 
symbol, and assume that $a|_{[ \delta,1[} $ is  differentiable for some 
$\delta \in ]0,1[$,  $a'$ satisfies \eqref{1.1a} and, for some 
constant $C>0$, there holds the  bound 
\bea
|a'(r)| \leq \frac{C}{(1-r)\big( \log(1-r) \big)^2} \ \ \ \mbox{for} \
r \in ]\delta ,1[.   \label{1.3r}
\eea
Then,   $T_a$ is a  bounded operator $  
H^{\infty}_v \rightarrow   H^{\infty}_v$. Moreover, if 
\bea
\lim_{r \to 1} |a'(r)| (1-r)\big( \log(1-r) \big)^2 = 0    \label{1.3t}
\eea
holds, then $T_a$ is compact,  if  and only if $\lim_{r \to 1} a(r) = 0$.
\ENT

{\bf Proof.} We can assume that $a$ is real-valued (otherwise consider Re $a$ 
and Im $a$ separately). Assume \eqref{1.3r} holds. For all $r \in ]\delta ,1[$ 
we get by the change of the integration variable $ \log(1- s) =: x$ and  
 $dx / ds  = - 1  / (1-s)$ that 
\bea
\int\limits_r^1 |a'(s)| ds  \leq C \int\limits_r^1
\frac{1}{(1-s)\big( \log(1-s) \big)^2} ds =
C \!\!\!\!\! \int\limits_{- \infty}^{\log (1-r)}  \!\!\!\!\!
\frac{1}{x^2} dx = \frac{C}{|\log (1-r)| }.   \label{1.3s}
\eea
Thus,  we can extend $a$ as a continuous function
to $]\delta ,1]$ by defining 
\beas
a(1) = \int\limits_\delta^1 a'(s) ds + a (\delta) \ \big(
= \lim_{r \to 1} a(r) \ \big),
\eeas
and by \eqref{1.3s} we obtain for all $r \in ]\delta ,1[$
\bea
|a(r) - a(1)| = \Big| \int\limits_r^1 a'(s) ds \Big| \leq 
\frac{C}{|\log (1-r)| } .  \label{1.3u}
\eea
This means, the function $a-a(1)$ satisfies \eqref{1.2} so that
the Toeplitz operator $T_{a-a(1)}$  is bounded. Since $T_{a(1)}$ is a 
multiple of the identity, $T_a = T_{a-a(1)}+T_{a(1)}$ is bounded. 

If \eqref{1.3t} holds, then we can repeat the calculation 
\eqref{1.3s}--\eqref{1.3u} so that the constant $C$ is replaced by a 
positive function $C(r)$ with $C(r) \to 0 $ as $r \to 1$. Then, we see from 
the analogue of \eqref{1.3u} that 
the function $a-a(1)$ even satisfies \eqref{1.3};
hence the operator $T_{a-a(1)}$  is compact, and if in  addition
$a(1)= 0$ then also $T_a$ is compact. If $\lim_{r \to 1} a(r)
= a(1) \not= 0$, then $T_a$ is a compact perturbation of a non-zero multiple of the 
identity which is not compact. \ \ $\Box$

\bigskip

All examples presented after Proposition \ref{prop1.3} also 
satisfy the assumptions of Theorem \ref{cor1.4}. 

The sufficient condition for the boundedness can be put into the following, very 
simple form. This  should be compared with corresponding results for non-radial 
symbols in [3]: we proved that for holomorphic $f$ on 
$\mathbb{D}$, the operator  $T_f$  is bounded, if and only $f$ is bounded while 
there are bounded harmonic  $g$ on $\mathbb{D}$ where $T_g$ is unbounded on $  
H^{\infty}_v$.

\BEC
If the symbol $a$ is radial and continuously differentiable on $[0,1]$,
then $T_a:   H^{\infty}_v \to   H^{\infty}_v$ is bounded. 
\ENC

For normal weights we can relax the assumptions on $a$ of Theorem \ref{th1.2} considerably.
\BET \label{th1.5}
Let $v$ be a normal weight. If $a \in L^1$
is radial and satisfies \eqref{1.2} then
$T_a$ is a bounded operator $  H^{\infty}_v \rightarrow  H^{\infty}_v$.

If $a$ satisfies \eqref{1.3} then $T_a$ is compact on $H_v^{\infty}$.
\ENT
We prove Theorem \ref{th1.5} in Section \ref{sec3}.
There is a variant of Theorem \ref{th1.2} for  exponential weights, too, 
without the restrictive smoothness requirements on $a$.

\BET  \label{th1.6}	
Let $v$ be an exponential weight of type $(\alpha,\beta)$. Assume that $a
\in L^1$ is radial and 
\bea
\label{1.4}
\limsup_{r \rightarrow  1} |a(r)|(1-r)^{- 1/2 - \beta/4} < \infty .
\eea
Then, $T_a$ is a bounded operator $  H^{\infty}_v \rightarrow  H^{\infty}_v$.
	
If
\bea  \label{1.5}
\limsup_{r \rightarrow  1} |a(r)|(1-r)^{- 1/2 -  \beta/4}
 = 0
\eea	
then $T_a$ is compact on $H_v^{\infty}$.
\ENT

We prove Theorem \ref{th1.6} in Section \ref{sec4}.

\section{Preliminaries.} 
\label{sec2}
To prove the theorems of Section \ref{sec1}  we need to recall some 
results of \cite{lu2} and \cite{toe}. We refer to these papers for a
more detailed exposition. 

Let $v$ be a weight on $\mathbb{D}$.  Fix $b > 2$. We define by induction  
the indices $0 \leq m_1 < m_2 < \ldots$ such that
\[  
b = \min \left( \Big( \frac{r_{m_n}}{r_{m_{n+1}}} \Big)^{m_n}
\frac{v(r_{m_n})}{v(r_{m_{n+1}})}, \Big( \frac{r_{m_{n+1}}}{r_{m_{n}}} \Big)^{m_{n+1}}
\frac{v(r_{m_{n+1}})}{v(r_{m_{n}})} \right). 
\]
This is always possible according to Lemma 5.1. of \cite{lu2}. (Actually it 
suffices to choose the indices such that the preceding minimum lies between 
$b$ and some constant $b_1 > b$.) Formula (6.1) of \cite{lu2} implies that 
\bea \label{2.1}  
\sup_n\frac{m_{n+1}-m_n}{m_n-m_{n-1}} < \infty 
\eea
so that we also have $\sup_n m_{n+1}/ m_n < \infty $ and 
$\sup_n (m_{n+1}-m_{n-1})/m_{n-1} < \infty .$

Now let  $h(\varphi) = \sum_{k \in \mathbb{Z}} b_k e^{ik \varphi}$
be a formal series with some numbers $b_k \in \bbC$ and $\varphi \in [0, 2\pi]$. Take the preceding numbers $m_k$  and define for every $n \in \bbN$ the operator
\begin{eqnarray}
(W_nh)(  \varphi) &=& \sum_{m_{n-1} < |k| \leq m_n} \frac{|k|-[m_{n-1}]}{[m_n]-[m_{n-1}]} b_ke^{ik \varphi}
+ \sum_{m_{n} < |k| \leq m_{n+1}} \frac{[m_{n+1}]-|k|}{[m_{n+1}]-[m_{n}]} b_ke^{ik \varphi} \nonumber \\
&=:&\sum_{k \in \mathbb{Z}} \beta_k  b_k e^{ik \varphi}  \label{2.0}
\end{eqnarray}
with coefficients $\beta_k = \beta_k(n)$ satisfying $|\beta_k| \leq 1 $ for
all $k$ and $n$.
Here $[r]$ is the largest integer not larger than $ r$. Obviously, $W_nh $ is 
always  a continuous function $[0, 2 \pi] \to \bbC$. The following
is Theorem 3.6. of \cite{toe}. 

\BET \label{th2.1}
Let the weight $v$ satisfy $(B)$. If $a \in L^1$ is radial
then $T_a$ is bounded as operator $H_v^{\infty} \rightarrow  
H_v^{\infty}$ if and only if 
\[ 
\sup_n \int_0^{2\pi}|(W_nf_a)(\varphi)| d \varphi < \infty
\]
and $T_a$ is a compact operator $H_v^{\infty} \rightarrow  
H_v^{\infty}$, if and only if 
\[ 
\lim\limits_{n \to \infty} 
\int\limits_0^{2\pi}|(W_nf_a)(\varphi)| d \varphi = 0. 
\]
Here, $f_a (\varphi)$ is  for  $\varphi \in [0, 2 \pi]$ the formal series
\[
f_a(\varphi) = \sum_{j =0}^{\infty}  \gamma_j e^{ij \varphi} \ \ \ \mbox{ with }  \ \ \gamma_n = \frac{\int_0^1r^{2n+1}v(r)a(r)dr }{\int_0^1r^{2n+1}v(r)dr}.
\]
\ENT
We recall that for radial symbols the Toeplitz operators reduces into 
a Taylor coefficient multiplier: if $h(z) = \sum_{n=0}^{\infty} h_n z^n$, 
then $T_a (z) = \sum_{n=0}^{\infty} \gamma_n h_n z^n  $. 

\bigskip

{\bf Examples.} If $v$ is normal, then one can
take $m_n = 2^{kn}$ for suitable fixed $k>0$ (see
\cite{lu2}, Example 2.4, and \cite{lul}).

For $v(r)= \exp(-\alpha/(1-r)^{\beta})$ one can take $m_n= 
\beta^2(\beta/\alpha)^{1/\beta}n^{2+2/\beta}-\beta^2 n^2$, and $r_{m_n}= 
1 -\left(\alpha/(\beta n^2)\right)^{1/\beta}$. This follows from (3.15), 
(3.16) and (3.30) of \cite{blt}. (There is a misprint in Theorem 3.1. of 
\cite{blt}, two times the exponent 2 is missing in the description of 
$m_n$.)	 

\BEC \label{cor2.2}
Let the weight satisfy $(B)$ and assume that $a \in L^1$ is radial and satisfies
$a|_{[s,1]}=0$ for some $s \in ]0,1[$. Then $T_a:H_v^{\infty} \rightarrow
H_v^{\infty}$ is compact.
\ENC

{\bf Proof.} We have
\begin{eqnarray*}
\left| \frac{\int_0^1a(r)r^kv(r)dr}{\int_0^1r^kv(r)dr} \right|  &   \leq &
\frac{\int_0^s|a(r)|r^kv(r)dr}{\int_{(1+s)/2}^1r^kv(r)dr} \\
               & \leq   &
 \left( \frac{2s}{1+s} \right)^k
 \frac{\int_0^s|a(r)|v(r)dr}{\int_{(1+s)/2}^1v(r)dr}.
\end{eqnarray*}
Hence, with $f_a$ as in Theorem 2.1,
\begin{eqnarray*}
\int_0^{2\pi}|(W_nf_a)(\varphi)| d \varphi  & \leq  &
c_1(m_{n+1}-m_{n-1}) \left( \frac{2s}{1+s}\right)^{m_{n-1}}  
\leq 
  c_2 m_{n-1}\left( \frac{2s}{1+s}\right)^{m_{n-1}} \end{eqnarray*}
for  universal constants $c_1, c_2$. Here we used (2.1). The right-hand side 
goes to $0$ as $n$ goes to $\infty$. Hence Theorem 2.1
finishes the proof. $\Box$

\bigskip

For $r > 0$ and an 
integrable function $f$ on $r \cdot \partial\mathbb{D}$ we put
\[ 
M_1(f,r) = \frac{1}{2 \pi}\int_0^{2 \pi} |f(r e^{i \varphi}) | d \varphi. 
\]
It is well-known that $M_1(f,r)$ is increasing with respect to $r$ if $f$ is a harmonic function.

Let $R$ be the Riesz projection, $ R: \sum_{k= - \infty}^{\infty} 
a_kr^{|k|}e^{ik \varphi} \mapsto \sum_{0}^{\infty} a_kr^{|k|}e^{ik 
\varphi}$. In the following we consider the Poisson kernel $p$, 
\[ 
p(re^{i\varphi}) = \sum_{k=-\infty}^{\infty} r^{|k|}e^{ik\varphi} \, ,
\ \ \ \mbox{where }re^{i\varphi} \in \bbD. 
\]
It is well-known that $p\geq 0$ and that
$M_1(p,r) =1$ for all $r \in [0,1[$. The following lemma will be needed 
later. 

\BEL \label{lem2.3}
Let $v$ satisfy condition (B) and consider the preceding numbers $m_n$ 
and operators $W_n$. Then we have
\[ 
\sup_n \sup_{0 \leq r < 1}M_1(RW_np,r) < \infty. 
\]
\ENL

{\bf Proof.}   According to Lemma 3.3  of \cite{lu2} we have
\begin{eqnarray*}
M_1(RW_np,r)  
& \leq & 4 \Big(\frac{[m_{n+1}]-[m_{n-1}]}{[m_{n}]-[m_{n-1}]} \Big)
\Big(3+4 \frac{[m_{n+1}]-[m_{n-1}]}{[m_{n+1}]-[m_{n}]}\Big)\\
& & 
\cdot\Big(1+\frac{[m_{n+1}]-[m_{n-1}]}{[m_{n-1}]} \Big) M_1(p,r). 
\end{eqnarray*}
Since $M_1(p,r) = 1$, the lemma follows in view of \eqref{2.1}. $\Box$

\section{Estimates for $\int_0^{2\pi}|(W_nf_a)(\varphi)| d \varphi $.}
\label{sec3}

For the proofs of the theorems of Section \ref{sec1} we will need the 
following estimate. 

\BEP \label{prop3.1}
Let $v$ be a weight which satisfies (B) and let $m_n$ be the numbers defined above
Theorem \ref{th2.1}. Assume that $a \in L^1$ is radial. Then there is a universal 
constant $c > 0$ with
\bea 
\label{3.1}
\int_0^{2\pi}|(W_nf_a)(\varphi)| d \varphi 
&\leq &c \log(m_{n}) \cdot \Bigg( \left|\frac{\int_0^1 a(r)
r^{2[m_{n-1}]+1}v(r)dr}{\int_0^1 r^{2[m_{n-1}]+1}v(r)dr} \right|
\rowpl
\sum_{[m_{n-1}] < k \leq [m_{n+1}]} \left|\frac{\int_0^1a(r) r^{2k+1}v(r)dr}{\int_0^1
	r^{2k+1}v(r)dr}-\frac{\int_0^1a(r) r^{2k-1}v(r)dr}{\int_0^1
	r^{2k-1}v(r)dr}  \right| \Bigg)	
\eea
and
\bea
\label{3.2}
\int_0^{2\pi}|(W_nf_a)(\varphi)| d \varphi \leq
 c\log(m_{n})  \frac{\int_0^1 |a(r)|
r^{2[m_{n-1}]+1}v(r)dr}{\int_0^1 r^{2[m_{n+1}]+1}v(r)dr}
\eea
for all $n$ large enough.
\ENP

To prove Proposition \ref{prop3.1} we need  a lemma. Given $m \in \bbN$, let 
$Q_m$ be the following projection acting on formal series (cf. \eqref{2.0}),
\[
Q_m\Big( \sum_{l=0}^{\infty} b_le^{il \varphi} \Big) 
= \sum_{l=0}^{m} b_le^{il \varphi}.
\]
It is well-known (see for example (2.7) in Ch. I, \cite{tor}) that 
\bea
\int_0^{2 \pi}\Big| Q_m \Big(\sum_{l=0}^{\infty} b_le^{il \varphi} \Big)\Big| 
d \varphi \leq 
d \log m \int_0^{2 \pi} \Big|
\sum_{l=0}^{\infty} b_le^{il \varphi}\Big| d \varphi  \label{3.2a}
\eea
where the coefficients $b_l$ for example form an $\ell^2$-sequence so that
the sum and the integral on the right converge and $d > 0$ is a  constant independent of $m$.

\BEL \label{lem3.2}
Let 	$f(e^{i \varphi}) = \sum_{k=0}^n b_k e^{ik \varphi}$ for
some $b_k \in \bbC$  and $n \in \bbN$, and let  $g(e^{i \varphi}) = 
\sum_{k=0}^n \alpha_{k} b_k e^{ik \varphi}$  for some coefficients 
$\alpha_k \in \bbC$. Then, 
\bea \label{3.3}
\int_0^{2 \pi} |g(e^{i\varphi})|d \varphi \leq
c  \log n \Big(|\alpha_0|+\sum_{k=1}^n|\alpha_k - \alpha_{k-1}| \Big) \int_0^{2 \pi} |f(e^{i\varphi})|d \varphi
\eea
where $c > 0$ is a  constant independent of $n$ and $f$.
\ENL

{\bf Proof.} We obtain, with $\beta_j = \alpha_{j}-\alpha_{j-1}$ for
$j=1, \ldots, n$ and $\beta_{0} = \alpha_{0}$,
\[ g(e^{i\varphi}) = \sum_{j=0}^n \beta_j \sum_{k=j}^n b_k e^{ik\varphi}. \]
Hence
\[\int_0^{2 \pi}|g(e^{i\varphi})|d \varphi \leq \left(|\alpha_0|+\sum_{k=1}^n|\alpha_k - \alpha_{k-1}| \right) \sup_j
\int_0^{2 \pi} \big| \big( ({\rm id} - Q_{j-1})f\big)(e^{i\varphi}) \big| d \varphi \]
from which we infer \eqref{3.3}, in view of \eqref{3.2a}. $\Box$

\bigskip

{\bf Proof of Proposition 3.1.}  Let  $f_a$ be again as in Theorem
\ref{th2.1}.  We have
\[
(W_nf_a)(e^{i \varphi})= \sum_{k=[m_{n-1}]}^{[m_{n+1}]} \frac{\int_0^1
a(r) r^{2k+1}v(r)dr}{\int_0^1 r^{2k+1}v(r)dr} \cdot \beta_k e^{ik \varphi}
\]
for certain $\beta_k$ with $|\beta_k| \leq 1$ (where $\beta_{[m_{n-1}]} 
=\beta_{[m_{n+1}]}=0$; see \eqref{2.0}). Now put $h(r e^{i\varphi}) = 
\sum_{j=0}^{[m_{n+1}]-[m_{n-1}]} \beta_{j+m_{n-1}}r^je^{ij \varphi}$ so that $h$ 
is a polynomial, hence  a holomorphic function. We obtain
\[
M_1(h,r) \leq M_1(h,1) = M_1(W_nRp,1) \ \ \ \mbox{ for } r \leq 1,
\]
where $p$ is the Poisson kernel.
Since $W_nRp$ is a polynomial we clearly find a radius $r(n) \in [0,1[$
such that
\bea  
\label{3.4}
M_1(W_nRp,1)\leq 2 M_1(W_nRp,r(n)) \ \ \ \mbox{ for all } n.
\eea
We use Lemma \ref{lem3.2} with  $f(e^{i\varphi})= h(re^{i\varphi})$
for fixed $r$, $b_j = \beta_{j+[m_{n-1}]}$ and
$$
\alpha_j=  \frac{\int_0^1a(s) s^{2(j+[m_{n-1}])+1}v(s)ds }{ 
\int_0^1 s^{2(j+[m_{n-1}])+1}v(s)ds }
$$ 
and obtain
\beas
& & \int_0^{2 \pi} |W_nf_a(e^{i \varphi})| d \varphi 
\rowleq 
c\log([m_{n+1}]-[m_{n-1}])\bigg( |\alpha_0|+\sum_{k=1}^{[m_{n+1}]-[m_{n-1}]}|\alpha_k - \alpha_{k-1}| \bigg)M_1(W_nRp,1).
\eeas
Then Lemma \ref{lem2.3} proves \eqref{3.1}.

To show \eqref{3.2} we see that
\begin{eqnarray*}
& & \int_0^{2 \pi} |W_nf_a(e^{i \varphi})| d \varphi
\\
&\leq &
\int_0^1 \int_0^{2 \pi}  \bigg|\sum_{k=[m_{n-1}]}^{[m_{n+1}]}
\frac{ r^{2k+1} a(r) v(r)}{\int_0^1 s^{2k+1}v(s)ds} \cdot \beta_k e^{ik \varphi} \bigg| d\varphi dr \\
&= &  \int_0^1   |a(r)| r^{2[m_{n-1}]+1}v(r)
\\
& & \ \  \  \cdot     \int_0^{2 \pi} \bigg|\sum_{k=[m_{n-1}]}^{[m_{n+1}]} 
\frac{1 }{\int_0^1 s^{2k+1}v(s)ds} \cdot r^{2(k-[m_{n-1}])}\beta_k e^{i(k-[m_{n-1}]) \varphi} \bigg| d\varphi  dr.
\end{eqnarray*}
Now put 
$$
\tilde{h}(r e^{i\varphi}) = \sum_{j=0}^{[m_{n+1}]-[m_{n-1}]}
\beta_{j+[m_{n-1}]}r^je^{ij \varphi}.
$$ 
Again we obtain
\[
M_1(\tilde{h},r) \leq M_1(\tilde{h},1) = M_1(W_nRp,1) \ \ \ \mbox{ for } r \leq 1.
\]

We use Lemma \ref{lem3.2} with  $f(e^{i\varphi})= \tilde{h}(re^{i\varphi})$
for fixed $r$, $b_j = \beta_{j+[m_{n-1}]}r^j$ and
$$
\alpha_j=\Big( \int_0^1 s^{2(j+[m_{n-1}])+1}v(s)ds  \Big)^{-1}.
$$ 
Then  $\alpha_j$ is increasing and
\[ 
|\alpha_0|+\sum_{k=1}^{[m_{n+1}]-[m_{n-1}]}|\alpha_k - \alpha_{k-1}| = |\alpha_{[m_{n+1}]-[m_{n-1}]} |
\]

The preceding estimate and \eqref{2.1},  \eqref{3.3},  \eqref{3.4}
yield   constants $c_1, c_2>0$ with
\begin{eqnarray*}
& & \int_0^{2 \pi} |W_nf_a(e^{i \varphi})| d \varphi
\\	& \leq &
2 \pi c_1 \log( [m_{n+1}]-[m_{n-1}])
\int_0^1  \frac{|a(r)| r^{2[m_{n-1}]+1}}{\int_0^1s^{2[m_{n+1}]+1}v(s)ds}
\cdot v(r) M_1(\tilde{h},r)dr
\\
& \leq &
2 \pi c_2 \log( [m_{n}]) \int_0^1  \frac{|a(r)| r^{2[m_{n-1}]+1}}{\int_0^1s^{2[m_{n+1}]+1} v(s)ds} \cdot v(r)
M_1(\tilde{h},1)dr.
\end{eqnarray*}
Now, Lemma \ref{lem2.3} also  shows  \eqref{3.2}. $\Box$

\bigskip

We recall the following

\bigskip

\BEL \label{lem3.3}
Let $v$ be normal. Then there is a universal constant $c > 0$ such that,
for any $k,m$ with $0 < k \leq m \leq 2k$, we have
\[
\frac{\int_0^1 r^k v(r) dr}{\int_0^1r^mv(r)dr }  \leq c.
\]
\ENL

{\bf Proof.} This is Lemma 4.5. of \cite{toe}.

\bigskip

\BEL \label{lem3.4}
For a function $a:[0,1] \rightarrow \mathbb{C}$, $ \epsilon > 0$  and $ \delta \in [0,1[$ there are constants $c_1,c_2 >0$ with
\begin{eqnarray*}
(a) \ \ \ \ \ \ c_1\sup_{n \geq 1/(1-\delta)}
\sup_{\delta \leq r \leq 1}|a(r)|r^n \log n & \leq &
 \sup_{\delta \leq r < 1} |a(r) \log(1-r)| \hspace{6cm} \\
  & \leq &
  c_2 \sup_{n \geq 1/(1-\delta)}
  \sup_{\delta \leq r \leq 1}|a(r)|r^n \log n \\
(b) \ \ \ \ \ \ \ \ \ \  c_1\sup_{n \geq 1/(1-\delta)}
\sup_{\delta \leq r \leq 1}|a(r)|r^n  n^{\epsilon} & \leq  &
\sup_{\delta \leq r < 1} |a(r)|/ (1-r)^{\epsilon} \hspace{6cm} \\
  & \leq &
  c_2  \sup_{n \geq 1/(1-\delta)}
  \sup_{\delta \leq r \leq 1}|a(r)|r^n  n^{\epsilon}  
\end{eqnarray*}
\ENL

{\bf Proof.} Put
$r = 1 - 1/n$, $n \geq 1/(1-\delta)$, and observe that $1/(1-1/n)^n$ is bounded.
 $\Box$

\bigskip
	
{\bf Proof of Theorem \ref{th1.5}}.  The inequalities \eqref{1.2} or \eqref{1.3}  
and Lemma 3.4 imply that there is $\delta \in [0,1[$ such that $\sup_{\delta \leq r < 1}|a(r) r^n | \leq c_0/ \log n$
for all $ n> 1$ and some constant $c_0$. Without loss of generality,  we may
assume that $\delta =0$, otherwise we take $a_1= a\cdot 1_{[\delta,1]}$ 
instead of $a$ and use the fact that $a = a_1 + a_2$ where
$a_2 = a \cdot 1_{[0, \delta]}$ yields the compact operator $T_{a_2}$.  We apply Proposition \ref{prop3.1}.
At first, we get
\[
\frac{\int_0^1| a(r)| r^{2m_{n-1}+1}v(r)dr}{\int_0^1 r^{2m_{n+1}+1}v(r)dr}  \leq \Big(\sup_{0 \leq r < 1}|a(r)|r^{m_{n-1}}\Big)\frac{\int_0^1 r^{m_{n-1}+1}v(r)dr}{\int_0^1 r^{2m_{n+1}+1}v(r)dr}
\]
According to  \eqref{2.1}  we have $\sup_k(m_{k+1}-m_k)/(m_k-m_{k-1})
< \infty$. This implies that $m_{n-1}+1 \leq m_{n+1}+1 \leq 2^q
(m_{n-1}+1)$ for some $q \in \mathbb{N}$ which is independent of $n$. If we apply  Lemma \ref{lem3.3}  $q$ times we see that
\[
\frac{\int_0^1s^{m_{n-1}+1}v(s) ds}{\int_0^1 s^{2m_{n+1}+1}v(s)ds}  \leq c^q
\]
where $c$ is the constant of Lemma \ref{lem3.3}. By  \eqref{3.2}
this shows, for some constant $c_1$
\begin{eqnarray*}
 \int_0^{2\pi}|(W_nf_a)(\varphi)| d\varphi  & \leq & c_1  \log(m_n)\left(\sup_{0 \leq r < 1}|a(r)|r^{m_{n-1}}\right)c^q \\
        & \leq & \frac{\log(m_n)}{\log(m_{n-1})} c^q
\end{eqnarray*}
where we used \eqref{1.2} and Lemma \ref{lem3.4}(a). With \eqref{2.1}
we see that $ \sup_n \int_0^{2\pi}|(W_nf_a)(\varphi)| d\varphi < \infty$. If \eqref{1.3} holds then the same estimate
yields 
$$
\lim_{n \rightarrow \infty}
\int_0^{2\pi}|(W_nf_a)(\varphi)| d\varphi =0.
$$ 
So  Theorem \ref{th1.5}
follows from Theorem  \ref{th2.1}. $\Box$

\bigskip

In order to prove Theorem \ref{th1.2} we need

\BEL \label{lem3.5}
Let $v$ be a weight on $ \mathbb{D}$ and
let $a: [0,1] \rightarrow \mathbb{R}_+$ be
continuous and non-increasing. Then
\[ \frac{\int_0^1 a(r)r^k v(r)dr}{\int_0^1 r^kv(r)dr} \geq
\frac{\int_0^1 a(r)r^{k+1} v(r)dr}{\int_0^1 r^{k+1}v(r)dr} \ \ \ \ \mbox{ for all } \ \  k=1,2, \ldots.\]
\ENL

{\bf Proof.} For $t \in [0,1]$ put
\begin{eqnarray*}
 F(t) & = & \left( \int_0^ta(r)r^kv(r)dr \right)\left( \int_0^tr^{k+1}v(r)dr \right), \\
 G(t) & = &\left( \int_0^ta(r)r^{k+1}v(r)dr \right)\left( \int_0^tr^{k}v(r)dr \right).  
 \end{eqnarray*}
Then $F$ and $G$ are differentiable and the
mean value theorem yields $s \in ]0,1[$ with
\[ \frac{F(1) -F(0)}{G(1) - G(0)} = \frac{ F'(s)}{G'(s)}\]
(Here we can assume that $a$ is not the zero function.)
Hence
\beas & & \frac{\left( \int_0^1a(r)r^kv(r)dr \right)\left( \int_0^1r^{k+1}v(r)dr \right)}{\left( \int_0^1a(r)r^{k+1}v(r)dr \right)\left( \int_0^1r^{k}v(r)dr \right)} 
\roweq
\frac{\left( \int_0^sa(r)r^kv(r)dr \right) s^{k+1}v(s) + a(s)s^kv(s)\left( \int_0^sr^{k+1}v(r)dr\right)}{\left( \int_0^sa(r)r^{k+1}v(r)dr \right) s^{k}v(s) + a(s)s^{k+1}v(s)\left( \int_0^sr^{k}v(r)dr\right)}.
\eeas
Since $a$ is non-increasing we have
\[
s^kv(s) \int_0^s\big( a(s)-a(r) \big) r^k(s-r)v(r) dr \leq 0
\]
which implies
\beas
& & \Big( \int_0^sa(r)r^kv(r)dr \Big) s^{k+1}v(s) + a(s)s^kv(s)
\Big( \int_0^sr^{k+1}v(r)dr\Big) 
\rowgeq
\Big( \int_0^sa(r)r^{k+1}v(r)dr \Big) s^{k}v(s) + a(s)s^{k+1}v(s)
\Big( \int_0^sr^{k}v(r)dr \Big). 
\eeas
Since $a$ is non-negative we obtain
\[ \frac{F(1) -F(0)}{G(1) - G(0)} \geq 1 \]
and hence 
\[\frac{\int_0^1 a(r)r^k v(r)dr}{\int_0^1 r^kv(r)dr} \geq
\frac{\int_0^1 a(r)r^{k+1} v(r)dr}{\int_0^1 r^{k+1}v(r)dr}. \hspace{4cm}\Box\]

\bigskip

{\bf Proof of Theorem 1.2.} Let the symbol $a$ satisfy the assumptions of Theorem 
\ref{th1.2}. Let us assume that 
\bea \label{3.5}
\limsup_{r \rightarrow 1} a'(r) < \infty  ,
\eea
otherwise we can consider $-a$. We may even assume that $a$ is differentiable on $[0,1[$.   Indeed put
\[ a_1(r) = \left\{\begin{array}{cc}
a(r)                            & r \in [\delta,1] \\
a(\delta) -a'(\delta)(\delta-r) & r \in [0,\delta]
\end{array} \right.\]
Then $a_1$ is differentiable on $[0,1[$.
Let $a_2= a \cdot 1_{[0, \delta]}$ and $a_3= a_1 \cdot 1_{[0, \delta]}$. According to Corollary \ref{cor2.2}, $T_{a_2}$ and $T_{a_3}$ are compact. Since $a = a_1 +a_2-a_3$ it suffices to assume (by perhaps taking $a_1$ instead of $a$) that $a$ is differentiable on $[0,1[$. 

Moreover it suffices to assume that $a$ is decreasing. Indeed otherwise consider
\bea 
\label{3.6}
\tilde{a}(r) = a(r) +d(1-r)
\eea
instead of $a$ where $ d> 0$ is a constant so large that $\tilde{a}' <0$ which 
exists in view of \eqref{3.5}. The symbol  $\tilde{a} $ satisfies 
\eqref{1.2} or \eqref{1.3}, too. If we have proved Theorem \ref{th1.2} for 
$\tilde{a}$ then it is also correct for
the symbol $1-r$. (Here take $a=0$ in \eqref{3.6}). So let us assume that
$a$ is differentiable everywhere, satisfies \eqref{1.2} or \eqref{1.3}
and is decreasing. Since  $\lim_{r \rightarrow \infty}a(r)=0$ we
obtain $a(r) \geq 0$ for all $r$.

We use again the terminology of Theorem \ref{th2.1}.
In view of Lemma \ref{lem3.5} with Proposition \ref{prop3.1}
we see that 
\bea
& & \int_0^{2\pi}|(W_nf_a)(\varphi)| d \varphi 
\rowleq 
c \log(m_{n}) \cdot \left(
 \frac{\int_0^1 a(r)
	r^{2[m_{n-1}]+1}v(r)dr}{\int_0^1 r^{2[m_{n-1}]+1}v(r)dr} +\frac{\int_0^1 a(r)
	r^{2[m_{n+1}]+1}v(r)dr}{\int_0^1 r^{2[m_{n+1}]+1}v(r)dr} \right) 
\label{3.7}
\eea
With \eqref{1.1} we obtain, for $k=m_{n-1}$ and 
$k=m_{n+1}$,
\[
\frac{\int_0^1 a(r)
r^{2k+1}v(r)dr}{\int_0^1 r^{2k+1}v(r)dr} \leq c_1 
\sup_r a(r) r^{k^{\epsilon}} 
\]
for a universal constant $c_1$. According to Lemma \ref{lem3.4}(a) and 
\eqref{1.2} we obtain a constant $c_2 > 0$  with 
$$ 
\sup_r a(r) r^{k^{\epsilon}} \leq 
 c_2 (\epsilon \log k)^{-1}.
$$ 
If we insert the last estimates in \eqref{3.7} we obtain $ \sup_n
\int_0^{2\pi}|(W_nf_a)(\varphi)| d \varphi
< \infty$ and we can apply Theorem \ref{th2.1}.
If we even have \eqref{1.3}, then the same estimates yield $\lim_{n 
\rightarrow \infty}  \int_0^{2\pi}|(W_nf_a)(\varphi)| d \varphi = 0$ and again 
Theorem \ref{th2.1} finishes the proof. \ \ $\Box$

\section{Exponential weights.} 
\label{sec4}
We now turn to the proof of Theorem \ref{th1.6}. Let us
fix an exponential weight $v$ of type $(\alpha,\beta)$, i.e. $v(r) = \exp(-\alpha/(1-r)^{\beta})$ for all $r$.
By analyzing the function $j \mapsto j^{\beta +1}-j^{\beta}$ we see that for every   $k > 0$ there is exactly one $j =j(k) > 1$ with
\bea \label{4.1} 
k = \alpha \beta (j^{\beta +1}-j^{\beta}). 
\eea
With this notation we see that $r_k :=1-1/j$ is the unique maximum point of the 
function $f(r) = r^k v(r)$.  Hence $f$ is increasing 
for $0 \leq r \leq r_k$ and decreasing for $ r_k \leq r < 1$. Moreover, in view of (4.1), there are constants $c_1, c_2 > 0$ with
\bea 
\label{4.2}  c_1 k^{1/(\beta +1)} \leq j \leq c_2 k^{1/(\beta +1)} \ \ \ \mbox{ for all } \ k \geq 1. 
\eea

\BEP
\label{prop4.1}
Let $k \geq1$ and the number $j =j(k)> 1$ be as chosen above and let $0 < \delta < 1$. Then 
there is a constant $d >0$, independent of $k$, such that
\[ 
\int_0^1 r^k \exp \Big(- \frac{\alpha}{(1-r)^{\beta}}\Big) dr \leq d
\int_{1-1/(\delta j)}^1 r^k \exp\Big(- \frac{\alpha}{(1-r)^{\beta}}\Big) dr. 
\]
\ENP

{\bf Proof.} Fix $1 >  \gamma > \delta $ and put
$x = 1-1/(\delta j)$, $y= 1- 1/(\gamma j)$. We may only consider large 
enough $j$ (and hence $k$) such that $0 < x$ and $\delta j > 1$.  Then we have $ 0 < x < y <1$. We use
\[ 
\exp\Big(- \frac{k}{t-1} \Big) \leq \Big(1-\frac{1}{t} \Big)^k \leq 
\exp\Big(-\frac{k}{t} \Big)
\]
whenever  $ 1 < t$. This implies
\begin{eqnarray*}
\int_0^x r^k \exp\Big(- \frac{\alpha}{(1-r)^{\beta}}\Big)dr
& \leq & x^{k+1}\exp\Big(- \frac{\alpha}{(1-x)^{\beta}}\Big) 
\\
& \leq & x^{k}\exp\Big(- \frac{\alpha}{(1-x)^{\beta}}\Big) 
\\
& \leq & \exp\Big(- \alpha \Big( \frac{\beta}{\delta}+ 
\delta^{\beta}\Big)j^{\beta} + \frac{\alpha \beta}{\delta}j^{
\beta -1}\Big)=: u_1.
\end{eqnarray*}
Moreover we have
\beas	
& & \int_y^{r_k} r^k \exp\Big(- \frac{\alpha}{(1-r)^{\beta}}\Big)dr
\geq  y^{k}\exp\Big(- \frac{\alpha}{(1-y)^{\beta}}\Big)\Big(1- 
\frac{1}{j} - y\Big) 
\rowgeq    
\exp\Big(- \alpha \beta \frac{j^{\beta+1}}{\gamma j -1}+ \alpha \beta 
\frac{j^{ \beta}}{\gamma j-1}  - \alpha  \gamma^{\beta}j^{\beta} \Big) 
\frac{1 -\gamma}{ \gamma j}  
\roweq   
\exp\bigg(- \alpha \Big( \frac{\beta}{\gamma}+ \gamma^{\beta}
\Big)j^{\beta} - 
\frac{\alpha \beta}{\gamma}j^{\beta -1} \Big(\frac{1}{\gamma 
-1/j}\Big)  \\
& & \ \ \ + \alpha \beta \frac{j^{\beta-1}}{\gamma-1/j} - 
\log \Big( \frac{\gamma j}{1- \gamma} \Big) \bigg)=: u_2. 
\eeas
Put $g(t) = \beta/t + t^{\beta}$ for $t > 0$. We easily see that $g$ is decreasing for $0 < t <1$.
Moreover 
\[ 
\frac{u_1}{u_2} \leq \exp\Big( -\alpha \big(
g(\delta) - g( \gamma) \big) j^{\beta} + d_1 j^{\beta -1} + d_1' \log j  \Big) 
\]
for some universal constants $d_1$, $d_1'$. Since 
$g(\delta) - g(\gamma) > 0$ this 
implies $ \limsup_{r \rightarrow \infty} u_1/u_2 <
\infty$. Hence
$u_1 \leq d_2 u_2$ for some universal constant $d_2$.
We obtain
\begin{eqnarray*}
& & \int_0^1 r^k \exp\Big(- \frac{\alpha}{(1-r)^{\beta}}\Big)dr  
\\
&  = &  
\int_0^x r^k \exp\Big(- \frac{\alpha}{(1-r)^{\beta}}\Big)dr +
\int_{x}^1 r^k \exp\Big(- \frac{\alpha}{(1-r)^{\beta}}\Big)dr 
\\
&  \leq & 
d_2 \int_{y}^{r_k} r^k \exp\Big(- \frac{\alpha}{(1-r)^{\beta}}\Big)dr + \int_{x}^1 r^k \exp\Big(- \frac{\alpha}{(1-r)^{\beta}}\Big)dr 
\\
& \leq  & 
(1+d_2) \int_{x}^1 r^k \exp\Big(- \frac{\alpha}{(1-r)^{\beta}}
\Big) dr.
\end{eqnarray*}
We finally put $d =1+ d_2$. $\Box$   

\BEC \label{cor4.2}
There is a constant $c > 0$ such that
\[ 	
\int_0^1 r^{k-k^{1/(\beta +1)}} \exp\Big(- \frac{\alpha}{
(1-r)^{\beta}}\Big)dr \leq c
\int_0^1 r^k \exp\Big(- \frac{\alpha}{(1-r)^{\beta}}\Big)dr 
\]
whenever $k \geq 1$. 
\ENC

{\bf Proof.} It is enough to consider sufficiently large $k$. Let $l = l(k)$ be such that
\[ 
k-k^{1/(\beta +1)} = \alpha \beta (l^{\beta +1} - l^{\beta}).
\]
For $k \geq k_0$, $k_0$ sufficiently large, there 
is a constant $c_0> 0 $ such that 
\[ 
(k - k^{1/(\beta+1)})^{1/(\beta +1)} \geq c_0 k^{1/(\beta +1)}. 
\]
Taking into account \eqref{4.2} for $k-k^{1/(\beta +1)}$
instead of $k$ we find a constant $c_1 > 0$ with
\bea  \label{4.3} 
l \geq c_1 k^{1/(\beta +1)}  .
\eea
Let $ \delta = 1/2$ and apply Proposition \ref{prop4.1} for
$ k - k^{1/(\beta+1)}$ instead of $k$. Together with \eqref{4.3} this yields
\begin{eqnarray*}
& & \int_0^1 r^{k - k^{1/(\beta+1)}} 
\exp \Big(- \frac{\alpha}{(1-r)^{\beta}} \Big) dr 
\\
& \leq &  d \int_{1-2/l}^1 r^{k - k^{1/(\beta+1)}} 
\exp \Big(- \frac{\alpha}{ (1-r)^{\beta}} \Big) dr 
\\ 
& \leq  &
\frac{d}{(1-2/l)^{k^{1/(\beta +1)}}} \int_0^1 r^{k } 
\exp \Big(- \frac{\alpha}{(1-r)^{\beta}} \Big) dr \\
&  \leq  & 
\frac{d}{\big(1-2 c_1 k^{- 1/(\beta +1)} \big)^{k^{1/(\beta +1)}}}
\int_0^1 r^{k } \exp \Big(- \frac{\alpha}{(1-r)^{\beta}} \Big) dr. 
\end{eqnarray*}
In order to complete the proof it is enough to take  $c$ such that 
\[ 
\frac{d}{\big(1-2 c_1 k^{- 1/(\beta +1)}  \big)^{k^{1/(\beta +1)}}} \leq c. 
\ \ \ \Box 
\]

\bigskip

Corollary \ref{cor4.2}  
proves \eqref{1.1} (with $\epsilon = 1/(\beta+1)$) for exponential weights 
and thus completes the proof of Proposition \ref{prop1.3} .

\bigskip

{\bf Proof of Theorem \ref{th1.6}}. By possibly taking $a \cdot 1_{[\delta,1]}$ 
instead of $a$ for suitable $\delta$ we can assume without loss of 
generality, in view of Lemma 3.4, that
$\sup_{0 \leq r \leq 1}|a(r)|r^k \leq c_0 /k^{1/2+\beta/4}$ for all $k$ and some constant $c_0$.  To obtain the indices $m_n$ of  Theorem
\ref{th2.1} we use \eqref{4.1} with $j = ( \beta n^2/ \alpha)^{1/\beta}$ (see (3.30), (3.15) and (3.16) of \cite{blt}). Hence
\bea 
\label{4.4}
m_n = \frac{\beta^{2+1/\beta}}{\alpha^{1/\beta}} n^{2+2/\beta} - \beta^2 n^2.
\eea
Let  $f_a$ be again as in Theorem \ref{th2.1}. We need to show that
\[
\sup_n \int_0^1 |W_nf_a(e^{i \varphi})| d \varphi < \infty.
\]
 We have
\[
(W_nf_a)(e^{i \varphi})= \sum_{k=[m_{n-1}]}^{[m_{n+1}]} \frac{\int_0^1
a(r) r^{2k+1}v(r)dr}{\int_0^1 r^{2k+1}v(r)dr} \cdot \delta_k e^{ik
\varphi}
\]
for certain $\delta_k$ with $|\delta_k| \leq 1$.

The equality \eqref{4.4} implies that there is
a constant $c_1 > 0$ such that
\bea \label{4.5}
(2k+1)^{1/(\beta +1)} \geq c_1 n^{(2+ 2/\beta)/(\beta +1)} = c_1 n^{2/\beta}
\eea
for all $k \geq   m_{n-1}$. Moreover, by  an application of
the mean value theorem to the function $n \mapsto m_n$ of \eqref{4.4}
we may assume that 
\bea 
\label{4.5a}  2(m_{n+1} - m_{n-1}) +1  \leq c_2 n^{1+2/\beta}.
\eea

The remark in the beginning of the proof, Corollary \ref{cor4.2}, \eqref{4.5} 
and the assumption \eqref{1.4} on $a$ yield a constant $c_3 > 0$ such that, for 
$m_{n-1} \leq k < m_{n+1}$, we have
\begin{eqnarray*}
& & \left|\frac{\int_0^1 a(r) r^{2k+1}v(r)dr}{
\int_0^1 r^{2k+1}v(r)dr}\right| \\
&   \leq &
\Big( \sup_r |a(r)|r^{(2k+1)^{1/(\beta +1)}} \Big)
 \cdot \frac{\int_0^1  r^{2k+1-(2k+1)^{1/(\beta +1)}}v(r)dr}{\int_0^1 r^{2k+1}v(r)dr} \\
&  \leq  &
c_3 \frac{1}{(n^{2/\beta})^{1/2 + \beta/4}} =   c_3 \frac{1}{n^{1/\beta +1/2}}.
\end{eqnarray*}
This implies by \eqref{4.5a}
\begin{eqnarray*}
\int_0^1 |W_nf_a(e^{i \varphi})| d \varphi
 & \leq &
\left(\int_0^1 |W_nf_a(e^{i \varphi})|^2 d \varphi \right)^{1/2} \\
& =     &
\Bigg(\sum_{k=m_{n-1}}^{m_{n+1}} \left|\frac{\int_0^1 a(r) r^{2k+1}v(r)dr}{\int_0^1 r^{2k+1}v(r)dr}\right|^2| \delta_k|^2 \Bigg)^{1/2} \\
& \leq  &
\Big((2(m_{n+1}-m_{n-1})+1) \frac{c_3^2}{n^{1+2/\beta}}\Big)^{1/2} \\
& \leq & c_2^{1/2}c_3.
\end{eqnarray*}
Hence, $\sup_n \int_0^{2\pi}|(W_nf_a)(\varphi)| d \varphi < \infty$.  So
Theorem \ref{th2.1} concludes   the first part of Theorem \ref{th1.6}.
If \eqref{1.5} holds then the same estimates as above show that
\hfill \break
$\lim_{n \rightarrow \infty} \int_0^{2 \pi} |W_nf_a(e^{i \varphi})| d
\varphi  =0$. Again, with Theorem \ref{th2.1} we see that then $T_a$ is
compact. \ \ $ \Box $

\noindent \textbf{Authors' addresses:}%
\vspace{\baselineskip}%

\noindent
Jos\'e Bonet: Instituto Universitario de Matem\'{a}tica Pura y Aplicada IUMPA,
Universitat Polit\`{e}cnica de Val\`{e}ncia,  E-46071 Valencia, Spain

\noindent
email: jbonet@mat.upv.es \\

\noindent
Wolfgang Lusky: Institut f\"ur Mathematik, Universit\"at Paderborn, D-33098 Paderborn, Germany.

\noindent
email: lusky@uni-paderborn.de \\

\noindent
Jari Taskinen: Department of Mathematics and Statistics, P.O. Box 68,
University of Helsinki, 00014 Helsinki, Finland.

\noindent
email: jari.taskinen@helsinki.fi

\end{document}